\documentclass[]{interact}

\usepackage{epstopdf}

\usepackage[numbers,sort&compress]{natbib}
\bibpunct[, ]{[}{]}{,}{n}{,}{,}

\usepackage{tikz}
\usetikzlibrary{patterns}
\begin{document}

\title{The polycons: the sphericon (or tetracon) has found its family}

\author{
\name{David Hirsch\textsuperscript{a} and Katherine A. Seaton\textsuperscript{b}\thanks{CONTACT K. A. Seaton. Email: k.seaton@latrobe.edu.au}}
\affil{\textsuperscript{a} Nachalat Binyamin Arts and Crafts Fair, Tel Aviv, Israel;\\ \textsuperscript{b} Department of Mathematics and Statistics, La Trobe University VIC 3086, Australia}
}

\maketitle

\begin{abstract} This paper introduces a new family of solids, which we call \textit{polycons}, which generalise the sphericon in a natural way. The static properties of the polycons are derived, and their rolling behaviour is described and compared to that of other developable rollers such as the oloid and particular polysphericons. The paper concludes with a discussion of the polycons as stationary and kinetic works of art.

\end{abstract}


\section{Introduction}

In 1980 inventor David Hirsch, one of the authors of this paper, patented `a device for generating a meander motion' \cite{DH}, describing the object that is now known as the sphericon. This discovery was independent of that of woodturner Colin Roberts \cite{PR}, which came to public attention through the writings of Stewart \cite{St}, P\"{o}ppe \cite{PS} and Phillips \cite{TP} almost twenty years later. The object was named for how it rolls ---  overall in a line (like a sphere), but with turns about its vertices and developing its whole surface (like a cone). It was realised both by members of the woodturning  \cite{SS, Sp} and mathematical  \cite{AT2, RK} communities that the sphericon could be generalised to a series of objects, called sometimes polysphericons or, when precision is required and as will be elucidated in Section \ref{roll}, the $(N,k)$-icons. These objects are for the most part constructed from frusta of a number of cones of differing apex angle and height.

\begin{figure}
\begin{center}
\resizebox*{16cm}{!}
{\includegraphics{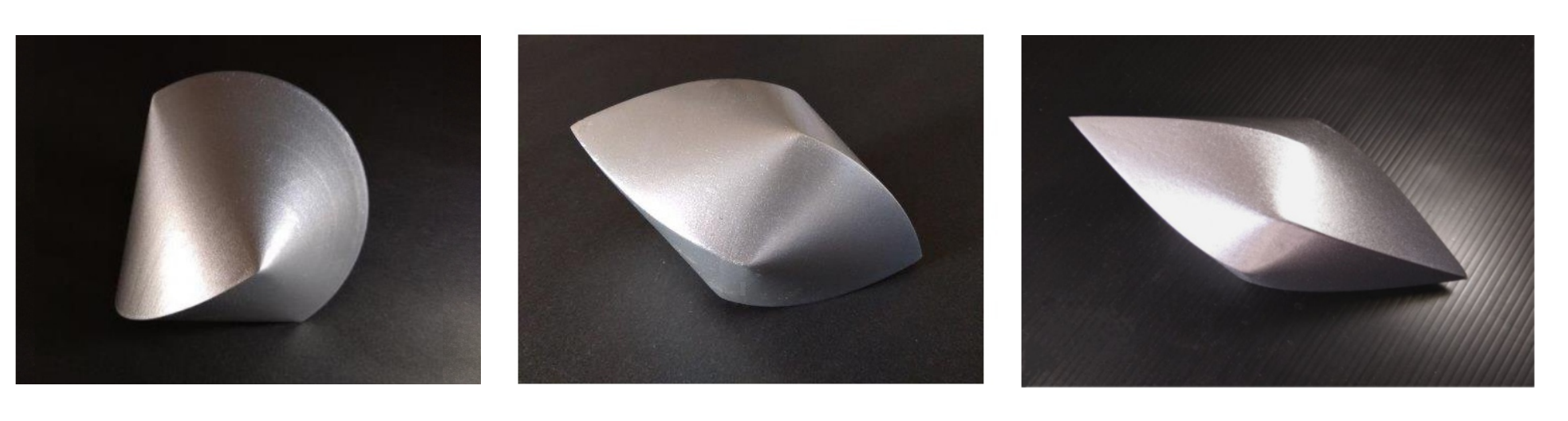}}
\end{center}
\caption{Models of a number of polycons, constructed by Hirsch. From left to right, they have $n=2$, $n=3$, $n=4$. Their sizes are, respectively, $10 \times 10 \times 10$ cm, $15\times 9 \times 7.5$ cm and $17\times 7 \times 6.5$ cm.} \label{photo}
\end{figure}

Now Hirsch has devised a different family, and one that is in many ways both more natural and more elegant, to which the sphericon gives birth. This series of generalisations we call the \emph{polycons}, each one being made from pieces cut from multiple identical cones. Figure \ref{photo} shows photos of working models of a number of polycons.

In this paper we present the theoretical construction of the polycons (Section \ref{cons}), and discuss their static mathematical properties (Section \ref{stat}). It is perhaps in their rolling that the polycons best display their elegance, and this motion is discussed in Section \ref{roll}. In this context, the features of the polycons are compared with those of the oloid \cite{PS1, PS2} (for which the sphericon is sometimes mistaken) and the polysphericons. We introduce the term {\em developable rollers} which encompasses all these rolling bodies. The paper closes with a discussion of the aesthetic motivation that led to the discovery of this new family of solids, and an artist's statement that shows the vertices, edges, surface and motion of the polycons in another light.

\section{Constructing the polycons}\label{cons}
There are four steps to create the polycon labelled by $n$, which we also call the $n$-con:
\begin{itemize}
\item[I.] The starting point is to take a right circular cone, sitting on the horizontal $(x,y)$-plane  with its symmetry axis aligned with the $z$-axis. Let the cone have radius $R$ and apex angle $\pi-\frac{\pi}{n}$, which is the internal angle of a regular $2n$-gon. The sloped surface of the cone makes angle $\frac{\pi}{2n}$ with the horizontal plane. Such a cone is shown in Figure \ref{cone}.

\item[II.] Next, two wedges are cut from the cone by two planes which intersect the $(x,y)$-plane along the $y$-axis, and which each make angle $\frac{\pi}{2}-\frac{\pi}{n}$ with the horizontal plane. These wedges are discarded. The $y=0$ cross-section of the piece that is kept is shown in Figure \ref{kite}. The cross-section is kite-shaped, with `top' angle $\pi-\frac{\pi}{n}$ and `bottom' angle $\frac{2\pi}{n}$. It is marked by the $z$-axis into two isosceles triangles. The whole solid piece cut from the cone has one curved surface which is part of the surface of the original cone (and includes its apex) and two flat surfaces which meet in a straight edge along the $y$-axis. We will classify the conic sections that form its other two edges below. 

\item[III.] Now $(n-1)$ more pieces, identical to the one created in steps I and II, are made.  The $n$ pieces are `glued' together using the flat surfaces to create a spindle-shaped object, which has a regular $2n$-gon as its $y=0$ cross-section and has $D_{nh}$ symmetry. While this is not yet the polycon, it has the surface area and volume of the final object. It has $n$ curved edges which pass through alternating vertices of the $2n$-gon.

\item[IV.] The last step in the construction is to cut the object produced in step III into two halves using the $y=0$ plane, rotate one of the halves through $\frac{\pi}{n}$ and then re-glue the halves together. This gives the polycon; its $y=0$ cross section is still a regular $2n$-gon but now it has $D_{nd}$ symmetry. 
\end{itemize}
Animations of this construction process are available to view for $n=3$ (the hexacon), $n=4$ (the octacon) and $n=5$ (the decacon) in a short film exhibited at the Bridges Conference, Linz in 2019 \cite{DH2}.

\begin{figure}
\begin{center}
\begin{tikzpicture}[scale=0.8]
  \coordinate (A) at (-4,0);
  \coordinate (B) at (4,0);
  \coordinate (C) at (-3,-1);
  \coordinate (D) at (3,1);
\coordinate (E) at (0,2);
  \draw[->, densely dashed] (A) -- (B)node[ right]{$x$};
  \draw[->, densely dashed] (C) -- (D)node[right]{$y$};
\draw[->, densely dashed] (0,-1.5)--(0,2.5)node[ right]{$z$};
   \draw[semithick] (-3,0) --(E); 
\draw[semithick] (E) --(3,0);
\draw[semithick] (-3,0) arc (180:360:3cm and 0.5cm);
\draw[dashed] (3,0) arc (3:180:2.97cm and 0.4cm);
  \node at (3,.1)[above]{$R$};
\draw[semithick] (3,0)--(3,0.1);
\draw (-2.3,0) arc (0:33:0.7cm);
\node at (-2,0.3){$\tfrac{\pi}{2n}$};
\end{tikzpicture}
\end{center}
\caption{A cone as described in Step I of constructing the polycons. It has apex angle $\pi-\frac{\pi}{n}$.}\label{cone}
\end{figure}
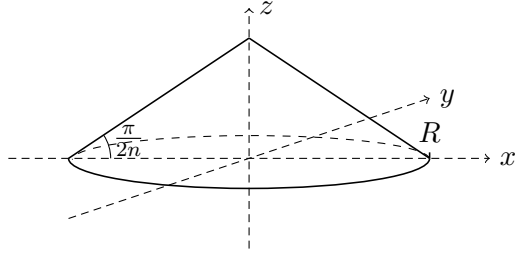
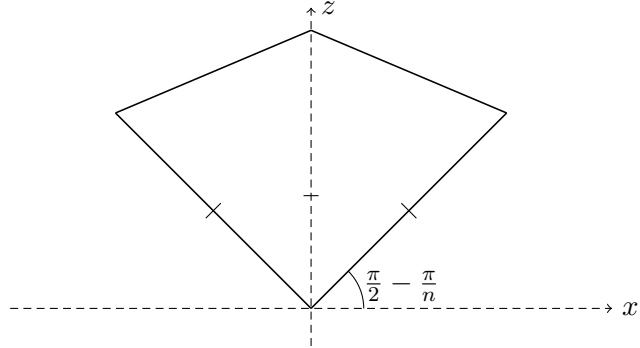
\begin{figure}
\begin{center}
\begin{tikzpicture}
  \coordinate (A) at (-4,0);
  \coordinate (B) at (4,0);
\coordinate (E) at (0,3.7);
  \draw[->, densely dashed] (A) -- (B)node[ right]{$x$};
\draw[->, densely dashed] (0,-.5)--(0,4)node[ right]{$z$};
  
   \draw[semithick] (-2.6,2.6) --(E); 
\draw[semithick] (E) --(2.6,2.6);
\draw[semithick] (0,0)--(-2.6,2.6);
\draw[semithick](0,0)--(2.6,2.6);
\draw (-1.4,1.2) --(-1.2,1.4);
\draw(1.2,1.4)--(1.4,1.2);
\draw(-.1,1.5)--(.1,1.5);
\draw (0.7,0) arc (0:45:0.7cm);
\node at (1.2,0.3){$\frac{\pi}{2}-\tfrac{\pi}{n}$};
\end{tikzpicture}
\end{center}
\caption{The $y=0$ cross-section of the piece kept from the cone in Step II is kite-shaped.}\label{kite}
\end{figure}
\subsection{The sphericon as a polycon}
When $n=2$, the four steps above describe an alternative construction of the sphericon, which is usually described as beginning with the creation of a bicone by rotating a square about its diagonal (see, for example, \cite{PR, St}). In step II above, the piece of the cone required is the whole  cone, and nothing is discarded, because the cutting planes coincide with the horizontal plane on which it sits. The conic section edges are two semicircles which make up the circular edge of the cone. The solid object formed in step III by glueing two cones together using their flat surfaces {\em is} the bicone, and step IV is the familar cut-rotate-reglue step \cite{KS} used in creating $(N,k)$-icons.  The sphericon is thus both the $2$-con and the $(4,1)$-icon. The left-most image in Figure \ref{photo} shows a model of this object.

\subsection{The vertices and edges of the polycon}
In general, the $n$-con has $2n$ curved edges, each one being half of the conic section created where the cone surface intersects one of the two cutting planes. On each side of the $y=0$ cross-section, $n$ edges of the polycon run (from every second vertex of the $2n$-gon) to a vertex situated on the $y$-axis and distance $R$ from the origin.  The edges on one side are offset by angle $\frac{\pi}{n}$ from those on the other side. The polycon has $2n+2$ vertices --- the two vertices at $y=\pm R$, which we call type A, and the vertices of the $2n$-gon, which we call type B. 
  
The nature of the conic section depends on $n$ in a way which we now make explicit. 
When a plane making angle $\gamma$ with the horizontal cuts the surface of a right circular cone with vertex angle $2\delta$, the resulting conic section has eccentricity $e=\frac{\sin(\gamma)}{\cos(\delta)}$. Thus  for the polycon edges
\begin{equation}
e=\frac{\cos(\frac{\pi}{n})}{\sin(\frac{\pi}{2n})}. \label{ecc}
\end{equation}

For the sphericon, $n=2$, the eccentricity is zero, so that the conic section is a circle.  (This is why there do not appear to be 4 edges and 6 vertices in the left-most image in Figure \ref{photo} as claimed above --- in this `new' description of the construction the edges are quarter-circles which meet smoothly on the $y$ axis to form two semi-circles. The two vertices of type A are indistinct.)

For the hexacon, $n=3$, the eccentricity is one, making the conic section  parabolic. The hexacon is thus another special case of the $n$-cons. 

 In equation (\ref{ecc}), the numerator increases with $n$ and the denominator decreases, so that for $n>3$ the conic sections are hyperbolic, having $e>1$. 

\subsection{Naming the family members}
 We have already informally introduced a naming convention (hexacon, octacon, decacon) for the polycons. They have $2n$ curved edges and a $2n$-gonal cross section.  It thus seems natural to use the Greek prefix that corresponds to $2n$  in naming the $n$-con.  (Hirsch originally considered calling the hexacon the `parabolicon' because of the nature of its edges, but this naming would not extend to the rest of the $n$-icon family.)

But according to this naming convention, the sphericon is now the odd one out! Just as a person may acquire a new name when joining a new family, or finding their true family, the sphericon can now be known as the  \textit{tetracon}. (It is, of course, in good company --- the cube has the alternative and accurate name of hexahedron.)

\section{Static properties of the polycons}\label{stat}

\subsection{Polycon volume}
\begin{figure}
\begin{center}
\begin{tikzpicture}[scale=1.2]
  \coordinate (A) at (-4,0);
  \coordinate (B) at (4,0);
  \coordinate (C) at (-3,-1);
  \coordinate (D) at (4,1.33);
\coordinate (E) at (1.5,3);

 \draw[ thick] (E)--(0,1);
\draw[semithick] (E) -- (.51,2);
\draw[thick] (E) --(3.52,1.52);
\draw[thick, pattern= dots] (0,1) arc (260:450:3cm and .5cm)--cycle;
\draw[semithick, densely dashed] (E)--(1.5,1.5);

\end{tikzpicture}
\end{center}
\caption{A generalised cone, as used in the volume calculation of the polycons. Its base (shaded) is bounded by a conic section and a line segment. It has one triangular face and one curved face. The height used in the volume calculation is indicated by the dashed line.}
\label{gen}
\end{figure}
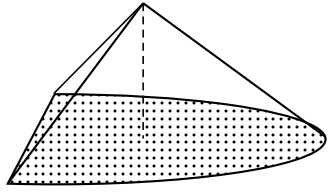

The object formed in Step III above has the same volume as the polycon. We can treat it as being comprised of $2n$ identical generalised cones as shown in Figure \ref{gen}. Each of them has as its apex  the apex of the cone it was cut from, and it has one flat surface (a triangle), and a curved developable surface, part of the surface of the original cone. The distance from the apex to the base of the generalised cone is $R \tan(\frac{\pi}{2n})\sin(\frac{\pi}{n})$. The base is part of the cutting plane bounded by the $y$-axis and a conic section.

As discussed by Apostol and Mnatsakanian \cite{AM}, it is  convenient to consider the projection of the conic section onto a horizontal plane, so that cylindrical polar coordinates can be used to calculate the area of the generalised cone base. The projected curves are of the same type (circle, parabola, hyperbola) with eccentricity $e^{\ast}$ given by
\begin{equation}
e^{\ast}=\tan({\delta})(\tan{\gamma})=\frac{\sin(\delta)}{\cos(\gamma)} e = \cot\left(\frac{\pi}{n}\right)\cot\left(\frac{\pi}{2n}\right). \label{ecc*}
\end{equation}
The volume of the polycon labelled with $n$ is given by
\begin{equation}
\text{Volume}=\frac{nR^3}{3}\tan\left(\frac{\pi}{2n}\right) I_n
\end{equation}
where
\begin{equation}
I_n=\int_{-\frac{\pi}{2}}^{\frac{\pi}{2}} \frac{1}{(1+e^{\ast}\cos(\theta))^2}d \theta. \label{int}
\end{equation}
These general expressions apply for all $n$, but the two smallest cases give simple answers (which can of course be found by other methods): \begin{equation}I_2=\pi \qquad \text{and} \qquad I_3=\frac{4}{3}. \label{simple}\end{equation}
The volume of the sphericon ($n=2$) is that of the bicone, $\dfrac{2\pi R^3}{3}$. \\For the hexacon ($n=3$) the volume is $\dfrac{4R^3}{3\sqrt{3}}$.

For $n>3$, so that $e^{\ast}>1$, the integral $I_n$ depends in a complicated manner on $n$:
\begin{multline}
I_n=\frac{(1-\cos(\frac{\pi}{n}))}{\cos(\frac{\pi}{2n}-\frac{\pi}{3})\cos(\frac{\pi}{2n}+\frac{\pi}{3})}\left[\frac{\cos(\tfrac{\pi}{n})}{2} + \phantom {\log_e \left(
\frac{1+2\sqrt{\rule{0ex}{1.8ex}\cos(\frac{\pi}{2n}-\frac{\pi}{3})\cos(\frac{\pi}{2n}+\frac{\pi}{3})}}{1-2\sqrt{\rule{0ex}{1.8ex}\cos(\frac{\pi}{2n}-\frac{\pi}{3})\cos(\frac{\pi}{2n}+\frac{\pi}{3})}}
\right)} \right.\\
\left. 
\frac{(1-\cos(\frac{\pi}{n}))^2}{8{\sqrt{\rule{0ex}{1.8ex}\cos(\frac{\pi}{2n}-\frac{\pi}{3})\cos(\frac{\pi}{2n}+\frac{\pi}{3})}}}\log_e \left(
\frac{1+2\sqrt{\rule{0ex}{1.8ex}\cos(\frac{\pi}{2n}-\frac{\pi}{3})\cos(\frac{\pi}{2n}+\frac{\pi}{3})}}{1-2\sqrt{\rule{0ex}{1.8ex}\cos(\frac{\pi}{2n}-\frac{\pi}{3})\cos(\frac{\pi}{2n}+\frac{\pi}{3})}}
\right)
\right]. \label{comp}
\end{multline}
This expression has been obtained by substituting (\ref{ecc*}) into the expression for $I_n$ (\ref{int}), which can be found in standard libraries of antiderivatives. Then by  application of trigonometric identities and some algebra, the expression has been rewritten as neatly as seems possible.

\subsection{Polycon surface area}
Using Theorem 10 of Apostol and Mnatsakanian \cite{AM} to calculate the  area of the curved surface of the object formed in step II above (part of the curved surface of the original cone), the total surface area of the polycon can also be expressed in terms of $I_n$: \[
\text{Surface area}= \frac{nR^2}{\cos(\frac{\pi}{2n})} I_n
\]
where $I_n$ is given in equations (\ref{simple}) and (\ref{comp}).
The surface area of the sphericon is $2\sqrt{2}\pi R^2$, and for the hexacon it is $\dfrac{8R^2}{\sqrt{3}}$.

\subsection{Inscribed solids} 

Consider the points on the conic section edges of the $n$-con at which they are intersected by the planes $y=\pm \frac{h}{2}$, where $0<h<2R$. Each plane intersects $n$ edges, giving the vertices of two parallel $n$-gons. The construction of the polycon ensures that the vertices of these two $n$-gons are offset by $\frac{\pi}{n}$, and they are `height' $h$ apart. Thus they are the vertices of an antiprism inscribed within the polycon; for appropriate choice of $h$ this antiprism will be uniform  (one in which the triangular faces connecting the $n$-gons are equilateral).

The relationship between the circumradius $b$ of the $n$-gon faces  and the height $h$ of a uniform antiprism is \cite{W}
\begin{equation}
h^2=\frac{b^2}{2}\left(\cos\left(\frac{\pi}{n}\right)-\cos\left(\frac{2\pi}{n}\right)\right). \label{height}
\end{equation}
Without loss of generality, we equate the polar coordinate expression for one of the conic section edges of the $n$-con with a particular vertex of the antiprism to give the condition
\begin{equation}
(\rho\cos(\theta), \rho \sin(\theta), \rho \cos(\theta) \cot(\theta))=(b\sin\left(\tfrac{\pi}{n}\right), \tfrac{h}{2}, b\cos\left(\tfrac{\pi}{n}\right)) \label{vertex}
\end{equation}
where, with $e^\ast$ given in (\ref{ecc*}) and $R$ being the radius of the cone used to construct the polycon,
\[
\rho=\frac{R}{1+e^\ast \cos(\theta)}.
\]
It is then possible to show that
\[
\cos(\theta)=\frac{2 \cos(\frac{\pi}{2n})}{\sqrt{3 +4 \cos(\frac{\pi}{n})}}
\]
and
\begin{equation}
b=\frac{R\sin(\frac{\pi}{2n})}{\sin^2(\frac{\pi}{2n})\sqrt{3+4\cos(\frac{\pi}{n})}+\cos(\frac{\pi}{n})\cos(\frac{\pi}{2n})}.
\end{equation}
The side-length of the $n$-gon is $a=2 b \sin(\frac{\pi}{n})$ and its height is given by (\ref{height}).

\begin{figure}
\begin{center}
\resizebox*{7cm}{!}{\includegraphics{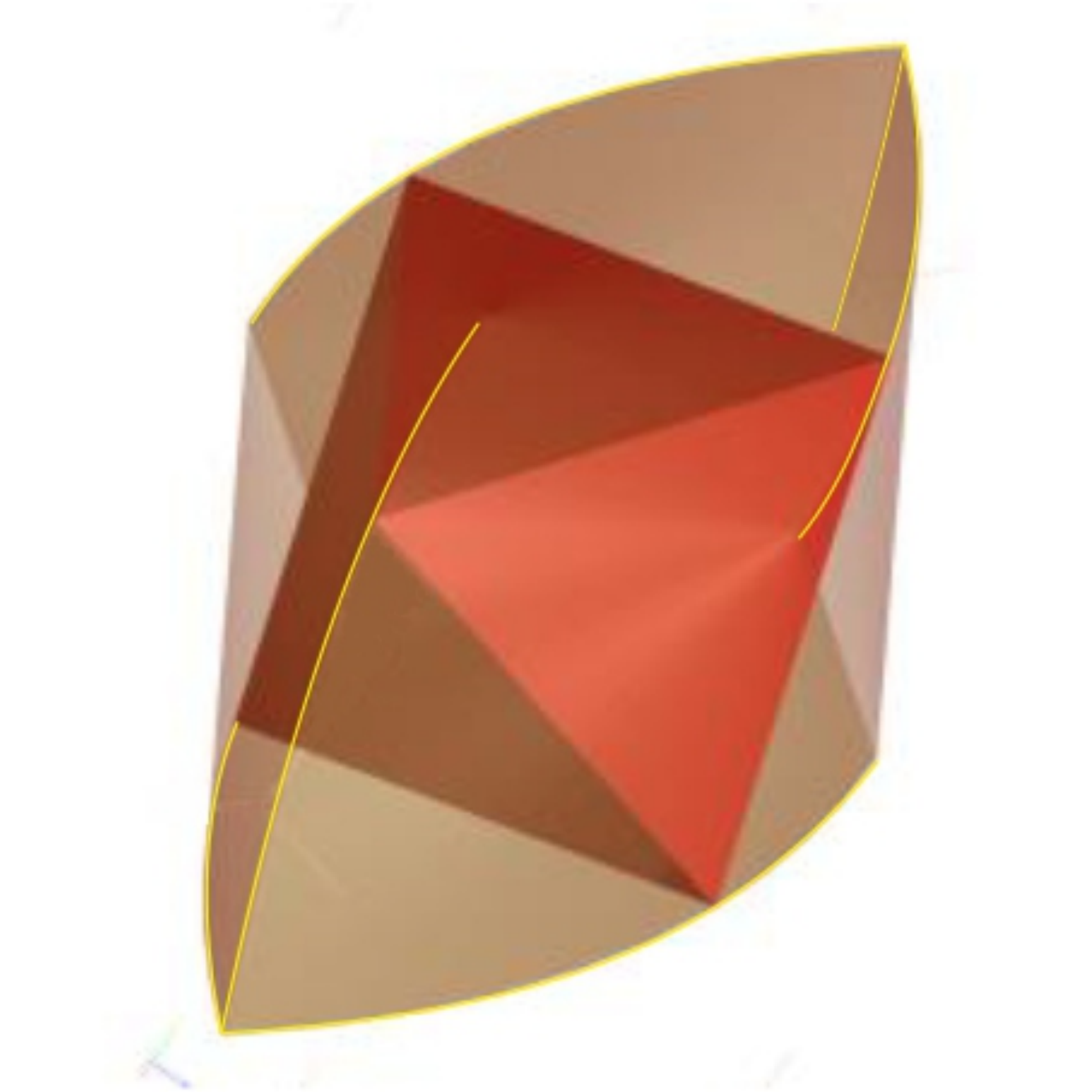}}
\end{center}
\caption{An octahedron (triangular antiprism) inscribed within a hexacon.}\label{inscribed}
\end{figure}

Again, $n=2$ and $n=3$ are special cases. Inside the sphericon is inscribed a degenerate antiprism, the tetrahedron (with side length $a=\frac{2\sqrt{2}R}{\sqrt{3}}$). Within the hexacon is a Platonic solid, the octahedron, with side length $a=R(\sqrt{15}-3)$. Figure \ref{inscribed} depicts this.

These are not the only solids inscribable within the $n$-cons, and it is intended to report on further aspects of this elsewhere. 

\section{Rolling properties of the polycons}\label{roll}
\begin{figure}
\resizebox*{7.5cm}{!}{\includegraphics{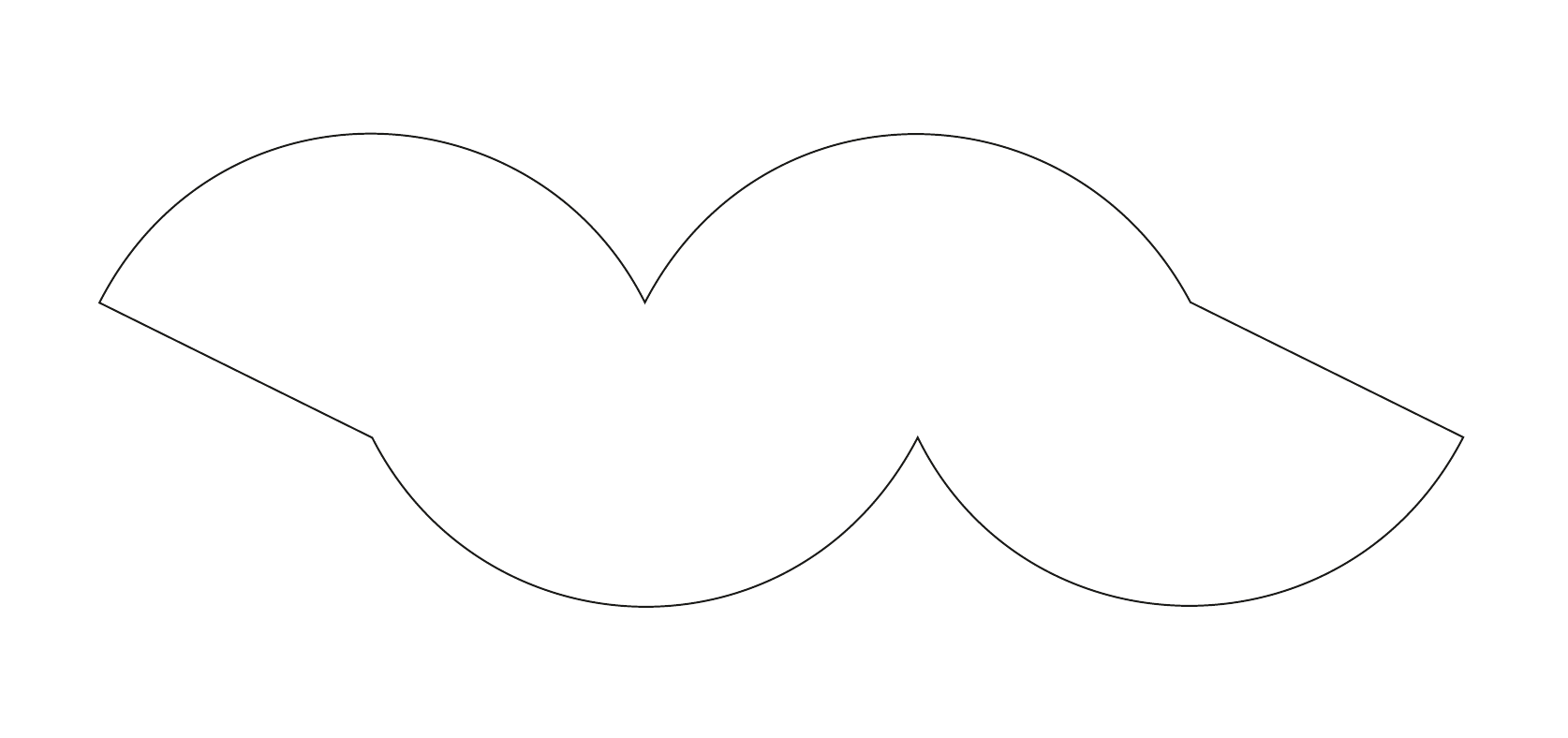}}\resizebox*{7.5cm}{!}{\includegraphics{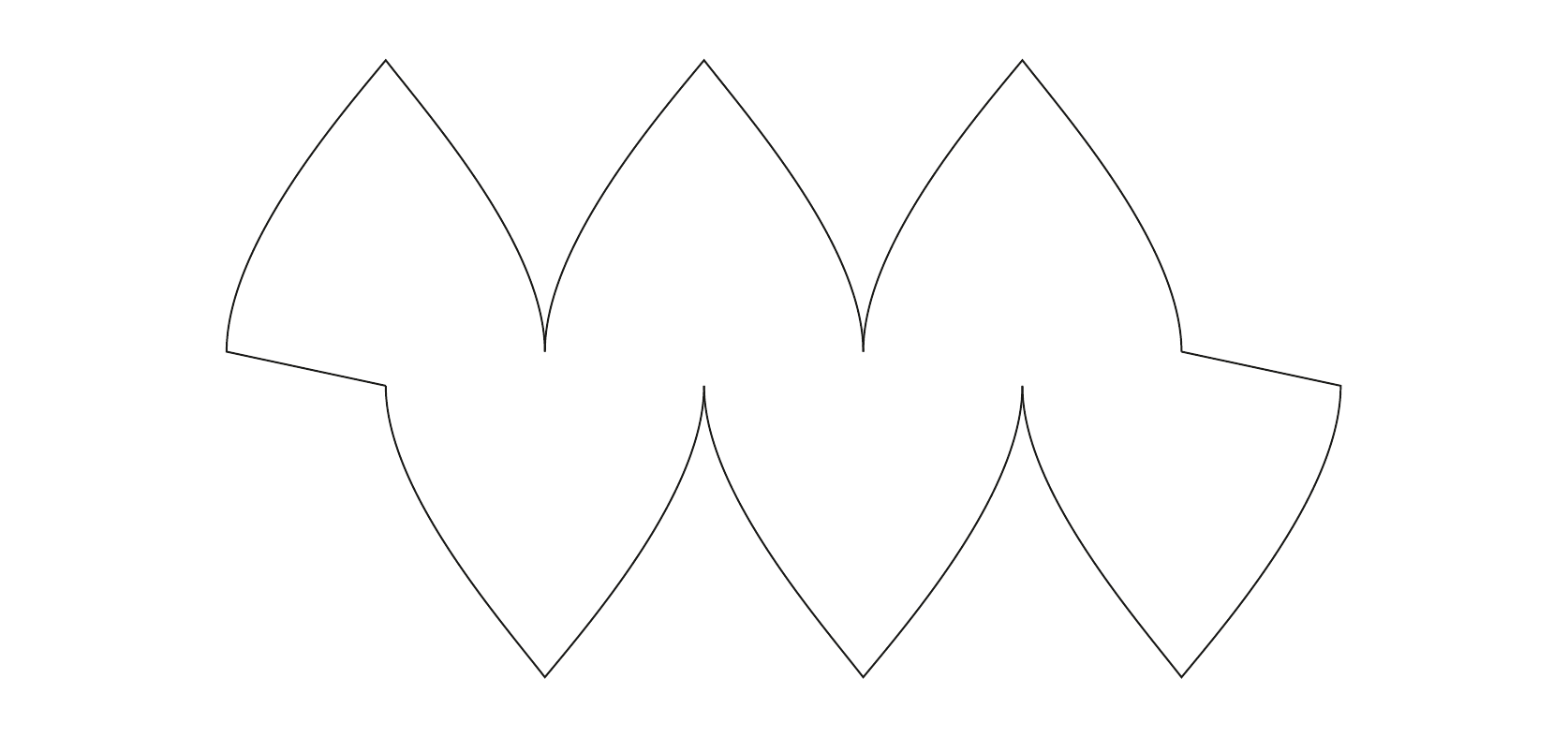}}
\caption{From left to right, the developed surface of the $n=2$  and $n=3$  polycon.} \label{mesh}
\end{figure}

The surface of each polycon is a single developable face. This face is shown (unrolled) in Figure \ref{mesh} for $n=2$ and $n=3$ and can be seen (unrolling) for the hexacon in the video \cite{DHf}. Thus the entire family has rolling properties that are related to the meander motion \cite{DH} or controlled wiggle \cite{St} of the sphericon, as do some members of the polysphericon ($(N,k)$-icon) family. 

We use the notation introduced in \cite{KS} for the $(N,k)$-icons: the index $N$ specifies an underlying $N$-gon, the index $k$ labels the twist used in construction and vv, mm and vm describe symmetry axes of the underlying polygon with v being `vertex' and m being `midpoint'. The members of the  $(N,k)$-icon family that have a single developable\footnote{\cite{RK, PR} use the word `continuous' and \cite{KS} uses `traceable'.} surface, and hence roll, are the vv $N$ even objects, when $k$  satisfies $\text{gcd}(\frac{N}{2},k)=1$.  These are also the particular polysphericons described as `prime' and identified as relevant to meander mazes in \cite{AT2}; \textit{prime polysphericons} will thus be used to specify this rolling subset of the $(N,k)$-icons. They consist of frusta of $\frac{N}{2}$ different cones. These cones come in pairs that match in terms of apex angle and height; when $N$ is an odd multiple of 2, the remaining section is actually from the cone with apex at infinity (i.e. a cylinder).  In forming the polysphericon, these various components are cut in half, twisted and rejoined. The cylindrical part of the surface has different rolling properties to the cones, and each pair of cones rolls differently from the others as well. Thus the rolling property of the prime polysphericons for $N>4$ is increasingly less fluid as $N$ increases --- something of a hiccupping wiggle with backtracking, rather than a meander.  Since the surfaces are developable, looking at templates for them gives some idea of this, as of course do videos of them moving. (See \cite{AT2, PS} for templates and, for example, \cite{pf,aero} for videos.)

In comparison to this `ugly duckling' motion, the polycons move like swans! 

Recall that the $n$-con consists of parts of $n$ identical cones, cut in half. The $2n$ vertices of type B are the apices of the original cones. Each edge of the $2n$-gon on the $y=0$ plane connects two type B vertices, being a segment of a generating line common to two cone pieces, one on each side of the $2n$-gon. When one of these edges makes contact with the flat surface, the rolling motion changes from being about one type B vertex to the next. As the $n$-con rolls on a surface, each of these vertices is in contact with the surface in turn for $\frac{1}{2n}$ of the time taken for a full revolution. The type A vertices (the only points from the edge of the original cone remaining once the polycon construction is complete) make contact with the surface instantaneously $n$ times each, alternately, during one revolution. Although the motion displays turns (about the type B vertices), because each cone section is identical the motion remains meandering, sinuous and smooth. It is best understood by watching it \cite{DH2}.

\begin{figure}
\begin{center}
\resizebox*{10cm}{!}{\includegraphics{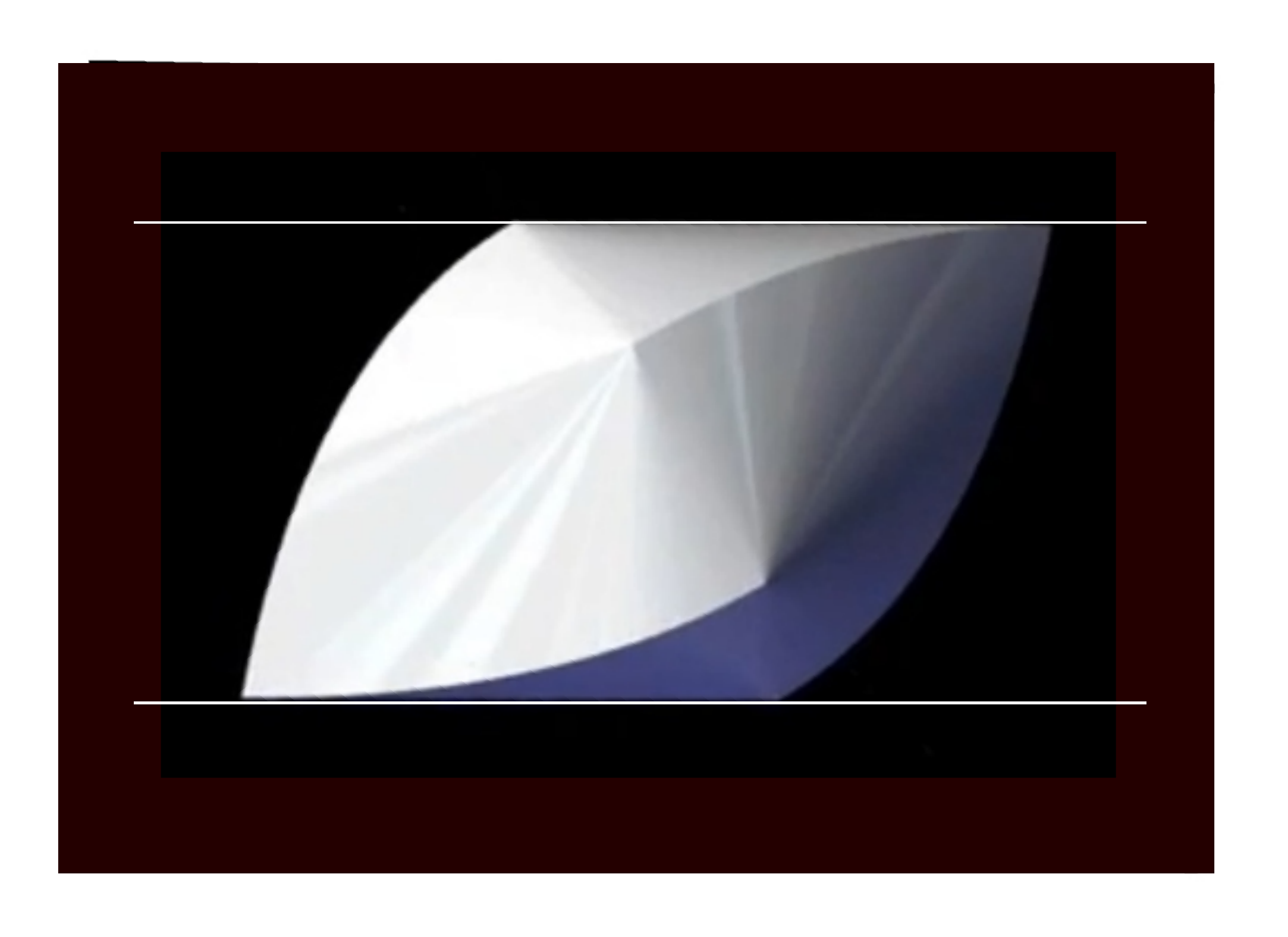}}
\end{center}
\caption{A computer-generated image of the hexacon, showing its constant height.} \label{hexa}
\end{figure}

The instantaneous line of contact between the polycon and the surface it is rolling on is a segment of one of the generating lines of a cone, and everywhere along this line the tangent plane to the cone (and hence the polycon) is the same. By symmetry, when $n$ is odd, this tangent plane is a constant distance from the tangent plane to the generating line on the polycon surface which is instantaneously uppermost (see Figure \ref{hexa}). Thus the polycons for $n$ odd are constant height rollers  (as is a right circular bicone, a cylinder or a prism with Reuleaux triangle cross-section). This feature is demonstrated in videos for $n=3$ and $n=5$ \cite{DHc}, and may be contrasted with $n=4$ \cite{DHnc}, which does not have this feature.

For all $n>2$, the length of the segment of the generating line in contact with the surface varies cyclically during the motion. This can also be seen in Figure \ref{hexa}, by considering the rays of light shown on the image. Its shortest length is the distance between two adjacent type B vertices: $2R\sec(\frac{\pi}{2n})\sin^2(\frac{\pi}{2n})$. Its longest length is the lateral side length of the original cone, the distance from a type A vertex to a type B vertex: $R\sec(\frac{\pi}{2n})$.  (These expressions also apply to the $n=2$ case (the sphericon), but the length of the generating line is constant: $\sqrt{2}R$.)

\subsection{Developable rollers}

As well as developing their entire surface when they roll (i.e., every part of it touches the surface on which it rolls), the prime polysphericons and the polycons share the property that the centre of mass remains at a constant height during rolling motion. The path  of the centre of mass (on the plane of constant height) is comprised of identical circular arcs for the polycons and of non-identical circular arcs (and, when $\frac{N}{2}$ is odd, line segments) for the polysphericons. 

\begin{figure}
\begin{center}
\resizebox*{10cm}{!}{\includegraphics{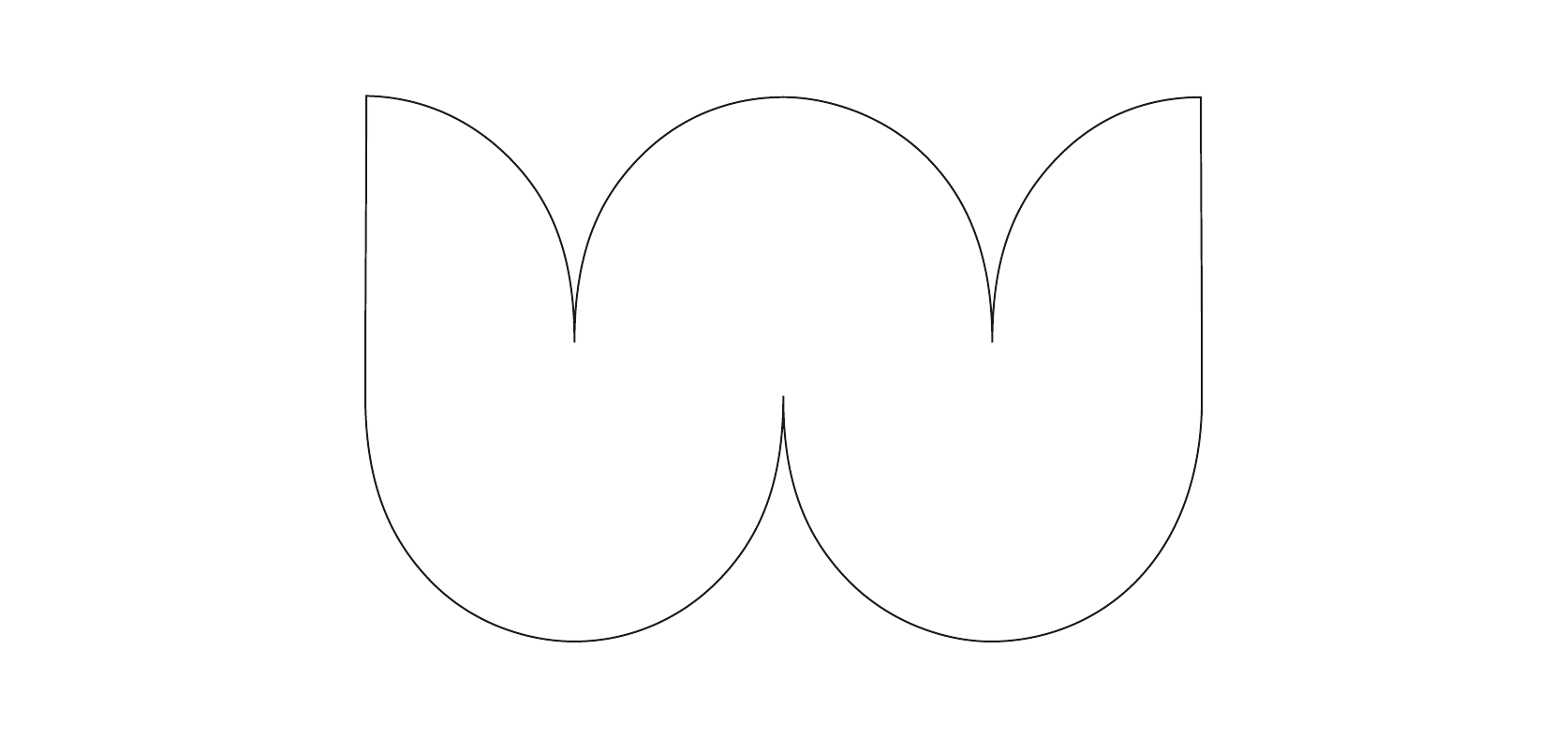}}
\end{center}
\caption{The developed surface of the oloid.} \label{oloid}
\end{figure}

The well-known oloid, described by Paul Schatz in 1929 \cite{PS1, PS2}, also has the property of developing its whole surface when it rolls, but its centre of mass does not remain at constant height \cite{DS}. The variation in its height is small, though, so that the rolling motion is pleasantly close to smooth.  The oloid is the convex hull of two congruent discs of radius $R$ set at right angles, with the centre of one lying on the edge of the other, i.e.,  the disc centres are distance $R$ apart. One third of each disc lies inside the oloid, so that it has no vertices, and it could be constructed from two sectors of discs with angle $\frac{4 \pi}{3}$. The developed surface of the oloid is shown in Figure \ref{oloid}.
Unlike the polycons in general, the line of contact between the oloid and the surface it rolls on (the generating line) is of constant length ($\sqrt{3} R$) \cite{DS}. The origin of its name has been suggested as being an abbreviation by Schatz of \textit{polysomatoloid}; alternatively, it has been suggested that its origin lies in the Greek word \textit{olos} (all)  \cite{Ferr}.

A summary of the properties of the objects discussed thus far in this section is given in Table 1.
\begin{table}\caption{Summary of  properties of the developable rollers: the $n$-cons (polycons), the prime polysphericons and the oloid}
\begin{tabular}{|l|c|c|c|}
\hline
Property&Polycons&Polysphericons&Oloid\\ \hline
Construction&Pieces of cong. & Frusta of&Convex hull of\\
& cones& non-cong. cones&two discs\\
Number of vertices&$2n+2$&$4$&$0$\\
Number of edges&$2n$&2&2\\
Nature of edges&conic sections& union of arcs& circular arcs\\ 
Height of centre of mass (rolling)&constant&constant &varies\\
Length of generating line segment&varies&fixed&fixed \\\hline
\end{tabular}
\end{table}

Hirsch's discovery of the (solid) sphericon \cite{DH} originated with two congruent half-discs, set at right angles with the centre point of their straight edges coincident. However, while this half-disc wobbler and its convex hull appear in \cite{EU, EU2},  the identification with the sphericon is not made (by name),  though its image appears and the `serpentine' motion of its centre of mass is depicted. The oloid, on the other hand, is discussed by name.  The oloid is thus another relative of the sphericon, related by way of the two-disc rollers. 

As an object comprising only two perpendicular circular discs of radius $R$, with centres separated by a distance $\sqrt{2}R$, the two-\textit{circle} roller was shown to roll with its centre of mass at constant height as long ago as 1966 \cite{ATS}. A mathematical description of its convex hull, paralleling exactly the discussion of \cite{DS} for the oloid, had to wait more than forty years. It appears in an unpublished and undated work \cite{HI}, which Ucke dates to 2011 \cite{EU2}. For completeness, we note that the convex hull of two arbitrary circles in 3 dimensions has been considered in \cite{NPSY}.

The [convex hulls of the] two-\textit{disc} rollers are generalisations which permit elliptical discs and half-discs \cite{EU, EU2}. With particular conditions on the ellipse dimensions and the separation of the disc centres (which the sphericon satisfies but the oloid does not), the centre of mass has constant height during rolling. 
The beautiful sculpture \textit{Rolling Lady} by Koman \cite{AKA}, which dates from the early 1980s, comprises two shallow bicones built on the two circles of such a roller, and which lie within the convex hull. While \textit{Rolling Lady} is constructed from developable surfaces and rolls, it is not what we term a developable roller.

\begin{figure}
\begin{center}
\resizebox*{8cm}{!}{\includegraphics{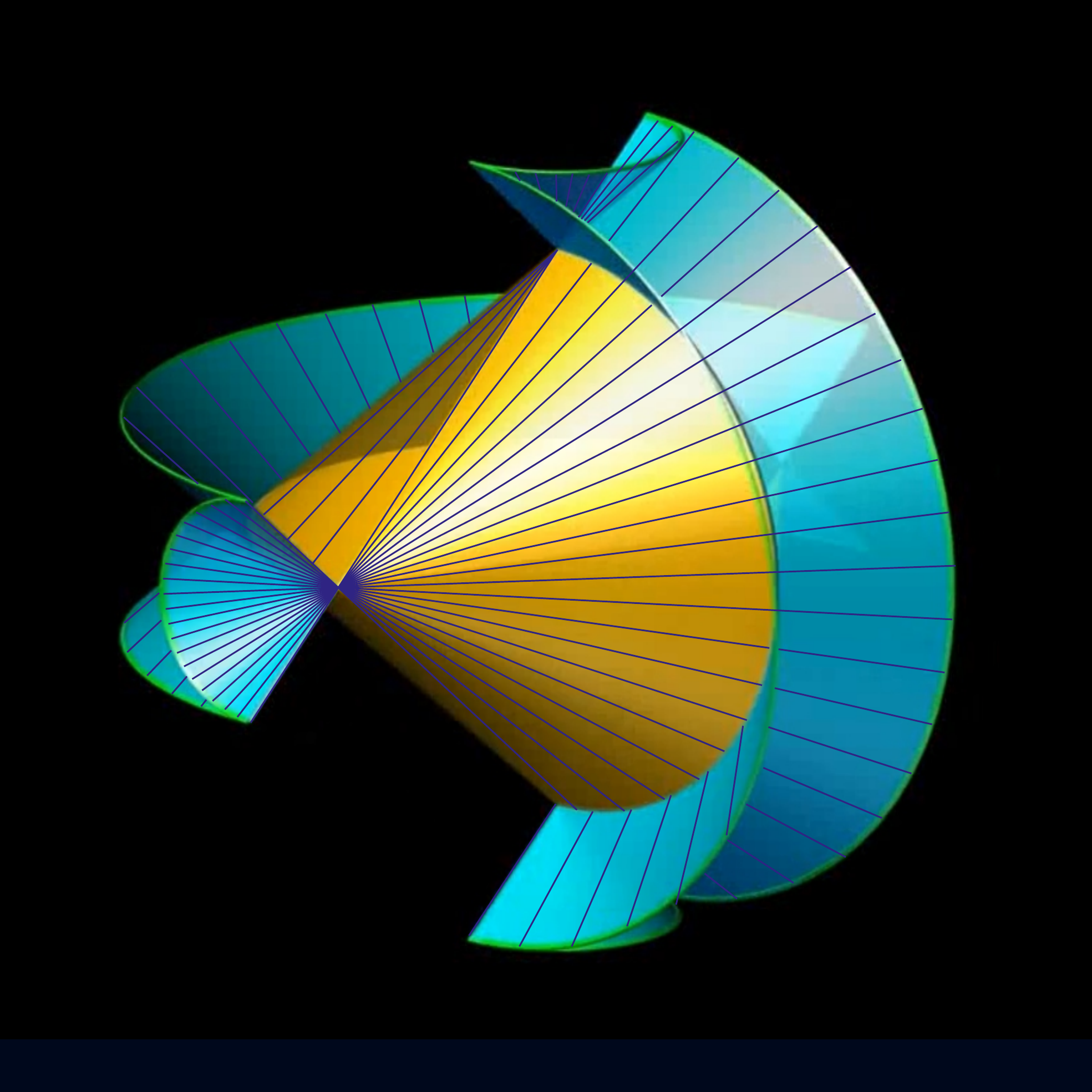}}
\end{center}
\caption{A ruled surface intersecting itself in such a way as to enclose a tetracon.} \label{roller}
\end{figure}

We define a developable roller as follows. Consider a ruled surface that intersects itself in such a way as to enclose a convex region of space, as shown in Figure \ref{roller} for a tetracon. Discard any part of the ruled surface that does not touch this closed volume. We call the body that remains a developable roller. Its surface is part of the original ruled surface and its edges are self-intersection curves of that ruled surface.  All the objects in Table 1, and the convex hulls of the two-disc rollers, are developable rollers. 

\section{The polycons as art}
\begin{figure}
\begin{center}
\resizebox*{15cm}{!}{\includegraphics{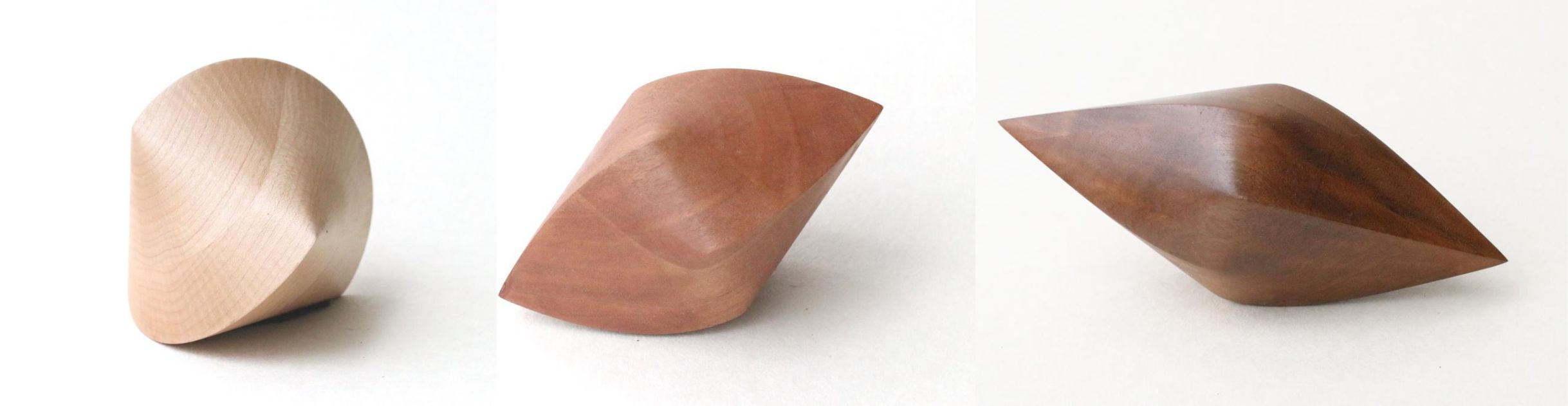}}
\end{center}
\caption{Polycons made for exhibition from maple, anigre and baillonella woods. Dimensions are $13 \times 13 \times 13$ cm,  $17\times10\times 8.5$ cm and $21\times 9 \times 8$ cm respectively.} \label{wood}
\end{figure}

The small working models shown in Figure \ref{photo} were made (in 2017) from aluminium using a CNC machine, hand-finished and coated with zinc paint. The larger wooden polycon sculptures of Figure \ref{wood}, which feature in the video \cite{DH2}, appear in the on-line gallery of the Art Show of the Bridges Conference, Linz  2019 \cite{DHbridges} for which they were created.
\subsection{The inspiration for the polycons}
It was not purely a mathematical motivation that led Hirsch to discover the polycons, but also the aesthetic and personal ones of an artist. In the words of an email he sent early in the collaboration that has led to this article, he wrote

\begin{quote}
After the publication of Ian Stewart's 1999 article, Steve Mathias published his series of sphericons that I believe you are familiar with. I liked this series, but there was something about it that bothered me. The surface of all the bodies in the series, except for the sphericon itself, is not uniform. It is always a combination of conical surfaces of various types, or of conical surfaces and cylindrical surfaces. I was busy with the thought, is it not possible to create a similar series in which the surface will be a uniform surface? .... It was clear to me that the surface I was looking for needed to be conical, and I realized that the polygon vertices would function as the apexes of the conical surfaces... From there the way was already short to complete the entire process...
And now to the name I have in mind [for the paper]: The polycons: the sphericon has found a family. 
As far as I'm concerned, the name has a double meaning. First of all, it is always interesting to examine something familiar from a broader perspective. There is a revelation of some truth, which was hidden and now it is revealed (or so I feel). The other aspect is personal. For me it is the completion of a work I started almost forty years ago... that somehow, reminds me of returning home to the family. 
\end{quote}

\subsection{The artist's statement}

The artistic qualities that can be attributed to bodies belonging to the polycon group are in fact inseparable from these qualities in the wider family of the developable rollers, of which the polycon group is a part. The uniqueness of the developable rollers is that they can be viewed as having two intertwined but separate aspects: their static physical appearance and the ``choreography" of their movement. 

Firstly, there are two features shared by all developable rollers: they all have a single face, and a structure that has rotational symmetry. What distinguishes the developable rollers from other single-faced bodies is that they have edges (and in most cases vertices too).

Single-faced bodies, like the sphere or the torus, have been associated throughout history with qualities such as perfection, unity, infinity and more. Another single-faced object, the M\"{o}bius ring, has become a symbol of recycling and renewal. To a large extent, all these values can also be attributed to the single-faced structure of the developable rollers. However, the fact that all developable rollers have edges doesn't allow them to be as ``perfect" as the sphere and the torus; rather, it makes their appearance much more complex, dynamic, interesting and alive.

It seems that there is no need to talk about the great importance of the use of symmetry in works of art from antiquity to our times. Beyond the innate attraction we have for symmetrical objects, symmetry is a symbol of the cosmic order present in our world. In contrast to  reflection symmetry that radiates a calm and static order, rotational symmetry expresses dynamism and motion. Indeed, the rotational symmetry of the developable rollers provides a strong sense of movement even when stationary.

In conclusion, we can say that the visual appearance of the developable rollers constitutes a unity of opposites. On one hand, there is the cyclic and continuous nature of their single face; on the other hand, the sharpness and dynamism of their edges, vertices, and rotational symmetry.

Secondly, in discussing the kinetic qualities of the movement of developable rollers, we can consider three main aspects: an element of surprise and complexity, a similarity to the movement of living beings, and a sense of rhythm.

Usually when we look at the movement of bodies that are not controlled by an intelligent factor, we are accustomed to being able to predict what the continuation of the motion will look like. This is the case when we roll a ball on the floor or throw a ball in the air. The characteristic meander or wiggle motion of the developable rollers may surprise someone who sees it for the first time. As their structure becomes more complex, so will their movement be more complex and unpredictable. This movement conveys a magical feeling of a structured ``plan" that seems to be implanted in their bodies.

The meandering movement, and especially the oscillation from side to side, gives rise to a strong association with the movement of living beings; perhaps to a walking duck or penguin or perhaps to someone tipsy, so that a comical feeling can also arise. The resemblance to living creatures makes the developable rollers an amusing thing that is easy to identify with.

When a ball or cylinder or even a cone rolls over a surface, except for a change of position in space, we cannot see any visual change. At every moment they look the same. With developable rollers,  this is different. At each moment of their rolling motion they reveal another part, different from the part that was revealed an instant ago. However, there is a cyclic pattern that repeats itself at intervals  characteristic to each specific developable roller. This is expressed more strongly as the structure of the developable rollers becomes more complex. This combination of constant change and repetition creates a sense of rhythm that is characteristic of dance, music and of course, kinetic art.

The visual uniqueness of the polycons presented in this paper, compared to other developable rollers, is their elongated structure and the sharp vertices at their two distant extremes. (In this respect, the first body in the series, the tetracon, is different from the rest.) Their shape is very similar to the trapezohedral shape, characteristic of many crystals. The pointed vertices also add to the sense of rhythm created during their rolling movement because they touch the rolling surface and move away from it in a regular cycle and thus serve as the ``rhythm section" of a musical band.

In conclusion developable rollers, and among them the bodies of the polycon group, have all the characteristics attributed to kinetic art. However, and this is also their uniqueness, their static appearance has qualities that allow them to serve as sculptural objects with an impressive presence even without moving.


\end{document}